\documentclass[12pt]{article}
\usepackage{amssymb, url, amsthm, amsmath}
\input cells.sty
\input rotate.tex

\usepackage{textcomp}
\cellsize=1em {
\def\boxit#1{
{\hbox{\lower3pt\hbox{\vrule\vbox{\hrule\kern2pt%
\hbox{\kern2pt$#1$\kern2pt}\kern2pt\hrule}\vrule}}}}
\bigskip
\bigskip
\hoffset =-0.7in \voffset =-0.5in \textheight =9in \textwidth =6.5in

\def\be{\begin{equation}}
\def\ee{\end{equation}}

\def\R{{\sf I\kern-.2em R}}
\def\N{{\sf I\kern-.2em N}}
\def\C{\kern.1em{\raise.47ex\hbox{$\scriptscriptstyle
$}}\kern-.40em{\sf C}}
\def\Z{{\sf Z\kern-.32em Z}}

\def\be{\begin{equation}}
\def\ee{\end{equation}}

\newtheorem{theorem}{\noindent Theorem}
\newtheorem{lemma}{\noindent Lemma}

\newtheorem{statement}{\noindent Proposition}

\urldef\vershik\url{vershik@pdmi.ras.ru}

\author {A.~M.~Vershik\thanks{%
St.~Petersburg Department of Steklov Institute of Mathematics.
E-mail: \vershik. Partially supported by the grants
RFBR 05-01-0089 and CRDF RUM1-2622-ST-04.%
}}

\title{A new approach to the representation theory of the symmetric groups, III: Induced representations
and the Frobenius--Young correspondence.}

\begin{document}
\maketitle

\rightline{\it To my friend Sasha Kirillov}

\begin{abstract}
We give a new (inductive) proof of the classical Frobenius--Young correspondence
between irreducible complex representations of the symmetric group
and Young diagrams, using the new approach, suggested in \cite{OV, VO}, to determining this
correspondence. We also give linear
relations between Kostka numbers that follow from the decomposition
of the restrictions of induced representations to the previous
symmetric subgroup. We consider a realization of representations
induced from Young subgroups in polylinear forms and describe  its relation
to Specht modules.
\end{abstract}

\section{Introduction}
In the classical representation theory of the symmetric groups,
the following theorem, which can be found in
all existing books on this subject, plays a key role.
This theorem goes back to the works by the pioneers of this theory,
Frobenius, Young, and Schur, and is sometimes called Young's rule
or the Frobenius--Young correspondence
(see \cite{F, Curt}).

{\it Let $\lambda\vdash n$ be a diagram with $n$ cells filled with
$n$ objects of arbitrary nature (for instance, the numbers
$1,2, \dots n$), and let
$H=H(\lambda)$ (respectively, $V=V(\lambda)$)
be the subgroup of the symmetric group
${\frak S}_n$ (the Young subgroup) consisting of all permutations of objects
inside the rows (respectively, columns) of this diagram. Consider the
representations
$Ind_{H(\lambda)}^{{\frak S}_n} \textbf{1}$ and $Ind_{V(\lambda)}^{{\frak S}_n}
\textbf {sgn}$
induced from the identity and sign representations of these subgroups,
respectively, to the whole symmetric group. Then their decompositions into
irreducible components contain exactly one common irreducible representation
$\pi_\lambda$, which has a simple multiplicity.}

This common representation $\pi_\lambda$ is determined up to
equivalence by the diagram, because it does not depend on by what
objects (or numbers) and how its cells are filled, so that different
fillings of the same diagram generate equivalent representations.
{\it It is this representation $\pi_\lambda$ that is assumed to be
associated with the diagram $\lambda$}; it is the ``principal''
component of both induced representations. For distinct diagrams
$\lambda$, the representations $\pi_{\lambda}$ are nonequivalent;
thus when $\lambda$ ranges over the set of all diagrams with $n$
cells, these representations exhaust the list of all classes of
nonequivalent irreducible complex representations of the group
${\frak S}_n$, because both the number of Young diagrams with $n$
cells and the number of classes are equal to $p(n)$, Euler's number
of partitions of an integer $n$.

This fact is a basis for further development of the theory. However,
its proof is not at all obvious. The traditional
proofs use combinatorial constructions that are far from
representation theory, so that they do not elucidate the matter
and thus cannot be extended to other Coxeter groups.
In another, less elementary, approach one obtains it from
the general theory of characters of the symmetric group,
which uses the techniques of the theory of symmetric functions.

Another disadvantage of such a method of defining this correspondence
is that it is implicit, and this fact predetermines the further steps
of the theory based on it. It is also worth observing
that the dimensions $\dim[Ind_{H}^{{\frak S}_n} \textbf{1}]$
and $\dim[Ind_V^{{\frak S}_n} \textbf {sgn}]$
of both induced representations
are much larger than the dimension $\dim\pi_{\lambda}$
of the main part of their intersection,
which is, perhaps, an evidence that any
proof of the fact under consideration
must be rather involved. The depth of this fact leads to
various combinatorial connections and parallels
(von Neumann's lemma, Gale--Ryser
theorem, etc.), links to the theory of partitions, the theory of generating
functions and symmetric functions, etc. This fact cannot be omitted
in any presentation of the theory.

In \cite{V, OV, VO}, another development of the whole theory was started.
In this approach, diagrams and tableaux appear quite naturally;
namely, standard Young tableaux are points of the spectrum of the
commutative
Gelfand--Tsetlin algebra; the set of all points of the spectrum, i.e., the set
of tableaux occurring in the same irreducible representation
$\pi$ of the group ${\frak S}_n$, is the set of tableaux corresponding to the
same diagram, and {\it it is this diagram that we associate with the
irreducible representation $\pi$}. Of course, the correspondence between
the irreducible representations and diagrams coincides with the
classical correspondence described above, but the new method is of
completely different nature. The inductive approach reveals
other important properties of this correspondence, Young's seminormal
and orthogonal forms become natural,
the bijection is almost obvious, and, most importantly,
the parametrization of
representations by diagrams is explicit and the branching rule
for the restriction of a representation to the previous subgroup
becomes obvious. Besides, this method without any modifications applies
to the Hecke algebras and, with some stipulations, to other series of
Coxeter groups.

With this construction of the representation theory of the symmetric
groups, we can completely omit the fact discussed above, it is not
necessary for further development of the theory; however, as
mentioned above, it is of independent importance, because it is a
clue to the role of representations induced from Young subgroups in
the representation theory. In this paper, which is mainly of
methodological character, we {\it prove this correspondence
(Theorem~1) by the inductive method and analyze the relation between
the representations induced from a Young subgroup for two successive
symmetric groups}. This allows us to derive simple and apparently
new recurrence relations between multiplicities of irreducible
representations, i.e., between Kostka numbers. In the last section,
we consider a realization of induced representations in spaces of
polylinear forms (tensors), describe a link to the classical Specht
modules, and give concrete examples. These questions, as well as a
number of new problems in combinatorics and the representation
theory of the symmetric groups, will be considered in detail
elsewhere.

I am grateful to the Schr\"odinger Institute (Vienna) and ETH
(Zurich) for invitations to give lecture courses on representation
theory in 2004 and 2005, and also to N.~V.~Tsilevich for
careful translation.

\section{Young--Frobenius correspondence}

We will proceed from the fundamental correspondence
$\lambda \leftrightarrow \pi_{\lambda}$, already obtained by the
inductive method (see \cite{OV,VO}):
$$
\{\mbox{Young diagrams with
$n$ cells}\} \Leftrightarrow \{\mbox{irreducible complex representations of
the group } {\frak S}_n\};
$$
In particular, it implies that the branching of representations of
the symmetric groups
${\frak S}_n$ is identical to the branching of Young diagrams,
and the branching graph is the graph of Young diagrams;
we will use this fact in what follows. The method of establishing this correspondence in
\cite{OV, VO} is based on considering not an individual symmetric group
${\frak S}_n$, but the whole inductive chain
$${\frak S}_1 \subset {\frak S}_2 \subset \dots \subset {\frak
S}_n,$$ uniquely determined up to isomorphism (for $n>6$), and on
the analysis of the Gelfand--Tsetlin algebra, naturally appearing in
this way, and its relations to the degenerate affine Hecke algebra.
Within this approach, the derivation of the Frobenius--Young
correspondence is simpler and more natural than in the conventional
presentation of the theory; it has the advantage that the branching
theorem precedes the more complicated theory of characters and so
on. The purpose of this paper is to explain, along the same lines,
the foundations of the classical method of constructing this
correspondence.

Before formulating the main theorems, let us recall the notation  and make
several preliminary remarks.

Let $\lambda \vdash n$ be an arbitrary Young diagram with $n$ cells. Let us
fix and denote by  $t_{\lambda}$ the tableau obtained by filling its cells
in an arbitrary way by the numbers
$1, \ldots, n$. All further considerations essentially depend only on the diagram,
the transition from a tableau to another tableau of the same shape
being equivalent to a conjugation in
${\frak S}_n$, so that our notation will involve only the diagram. The tableau
$t_{\lambda}$ determines two partitions of the set
$\{1, \ldots, n\}$: the partition into the rows of $t_\lambda$, denoted by
the same symbol
$$\lambda=(\lambda_1,\ldots ,\lambda_n),$$ and the partition
$$\lambda'=(\lambda'_1,\ldots, \lambda'_n)$$ into the columns of
$t_{\lambda}$. Denote by $H({t_\lambda})\equiv H(\lambda)$ the Young subgroup
associated with $t_\lambda$,
i.e., the subgroup of
${\frak S}_n$ consisting of permutations that preserve
the partition into the rows of ${t_\lambda}$:
$${\frak S}_{\lambda_1}\times \dots\times {\frak S}_{\lambda_n};$$
by $V(t_{\lambda})\equiv V(\lambda)$
we denote the Young subgroup preserving the columns of ${t_\lambda}$:
$${\frak S}_{\lambda'_1}\times \dots\times {\frak S}_{\lambda'_n}.$$
Obviously, $V(\lambda)=H(\lambda')$, where
$\lambda'$ is the diagram conjugate to $\lambda$. Given a group $G$,
by $Ind_H^G \pi$ we denote the operation of inducing
a representation $\pi$ of a subgroup $H$ to the group
$G$, and by $Res^G_K(\Pi)$ the restriction of a representation
$\Pi$ of the group $G$ to a subgroup $K$;
the one-dimensional identity and sign representations
of the symmetric group are denoted by
$\textbf{1}$
and $\textbf{sgn}$, respectively.

The notation $\lambda\succ \gamma $  or $\gamma \prec \lambda$ means
that a diagram $\lambda$  immediately follows
a diagram $\gamma$ (or $\gamma$ immediately precedes $\lambda$)
in the Young graph (i.e., $\lambda$ is obtained from $\gamma$
by adding one cell). The notation
$\lambda \unlhd \mu$ or $\mu\unrhd\lambda$ means that for every
$k$, the sum of the lengths of the first $k$ rows of
$\lambda$ does not exceed the sum of the lengths of the same rows of
$\mu$; this is the dominance ordering on diagrams with the same number of cells.

The main fact leading to the classical version of the correspondence
  $$\{\lambda\} \Leftrightarrow \{\pi_{\lambda}\}$$
between the diagrams with $n$ cells and the complex irreducible
representations of the symmetric group  ${\frak S}_n$ (the
Frobenius--Young correspondence) is as follows:
\newpage
\begin{theorem}[Young--Frobenius correspondence]
{$$Ind_{H(\lambda)}^{{\frak S}_n}[\textbf{1}] \bigcap
Ind_{V(\lambda)}^{{\frak S}_n}[\textbf{sgn}]=\{\pi_{\lambda}\}.$$}
\end{theorem}

The left-hand side of this formula should be understood as the intersection
of two multisets of irreducible representations that appear in the decomposition
of each of the induced representations associated with the diagram
(more exactly, with the partitions generated by this diagram,
i.e., the horizontal and vertical partitions into the rows and columns
of the diagram, respectively), taking into account
the multiplicities.
The claim is that the multiplicity of the unique irreducible
representation belonging to this intersection is equal to one
in each of the two representations, and it is this representation that
is associated with the diagram $\lambda$. In our approach,
the correspondence between irreducible representations and diagrams
is already established (by the above-mentioned method of
\cite{OV, VO} or in some other way),
and we must justify the classical version, i.e., give a simple
and conceptual proof that the above intersection  consists of a single
representation $\pi_{\lambda}$ which has multiplicity one.
We will prove this theorem by induction on the degree of the symmetric group.
For this, we will use well-known facts from the theory of induced representations,
which, for some reason, rarely appear in manuals on the  symmetric groups.
We start from the general and well-known Mackey's formula
for the restriction of an induced representation
to a subgroup.

\begin{theorem}
Let $G$ be a finite group,
$H$, $K$ be two subgroups of $G$, and
$\rho$ be a representation of
$H$ in a space
$W$. Let $H_s=sHs^{-1} \bigcap K$, and consider the representation
$\rho_s$ of the subgroup $H_s$ in the same space
$W$ defined by
$\rho_s(x)=\rho(sxs^{-1})$. Then
$$
Res_K Ind_H^G W =\sum_{s \in K\backslash G/H}Ind_{H_s}^K \rho_s.
$$
\end{theorem}

The sum in the right-hand side ranges over the space of double
cosets of the subgroups $H$ and $K$. It is not difficult to check
this formula directly by the definition of an induced representation
(see \cite{Se, Ma}).

The assertion of the following lemma is obtained by simply
applying Mackey's formula to our case:

\begin{lemma}
$$Res_{{\frak S}_{n-1}}Ind_{H(\lambda)}^{{\frak S}_n}
\textbf{1}=\sum_{\gamma:\gamma \prec \lambda} c(\lambda,\gamma) Ind_{H(\gamma)}^{{\frak
S}_{n-1}}\textbf{1}.$$ Here we denote by $c(\lambda,\gamma)$ the number of ways to obtain a
partition $\lambda \vdash n$ from partition $\gamma \vdash (n-1)$.
\end{lemma}

The coefficient $c(\mu,\lambda)$ can be defined in another way,
as the multiplicity of the row being modified in the diagram
$\mu$ (here diagrams are understood as partitions); or as the number
of ways to obtain a partition with diagram
$\gamma \vdash (n-1)$ as the restriction of a partition with diagram
$\lambda \vdash n$; or, finally, as the number of blocks in the partition
$\lambda$ such that decreasing one of them by one element yields the partition
$\gamma$.

\begin{proof}
To prove this lemma, let us apply Theorem~2 (Mackey's formula):
$$Res_K Ind^G_H\textbf{1}=\sum_{s\in K\backslash G/ H} Ind_{H_s}^K\textbf{1},$$
where  $H_s = sHs^{-1}\bigcap K$. In our case, $G={\frak S}_n$, $H=H(\lambda)$,
$K={\frak S}_{n-1}$, and $H_s$ ranges over all possible
Young subgroups in ${\frak S}_{n-1}$ obtained as the intersections of
${\frak S}_{n-1}$ with subgroups conjugate with $H(\lambda)$ in
${\frak S}_n$.

Thus $H_s$ ranges over the set of Young subgroups
$H({\gamma})$ of the group
${\frak S}_{n-1}$ that correspond to various partitions
$\gamma$ obtained by the restriction of
$\lambda$ to various subsets of cardinality
$n-1$ in the set $\{1,2,\ldots, n\}$. This gives the coefficient
$c(\lambda,\gamma)$ in the desired formula.
\end{proof}

In the proof of Theorem~1 we will use only part of information
contained in this lemma; namely, that every irreducible representation
$\pi_{\rho}$ of the group ${\frak S}_{n-1}$ appearing in
the restriction to ${\frak S}_{n-1}$
of the induced representation $Ind_{H(\lambda)}^{{\frak
S}_n}\textbf{1}$ appears in one of the induced representations
$Ind_{H(\gamma)}^{{\frak S}_{n-1}}\textbf{1}$, where
$\gamma \prec \lambda$.
We will return to the question of multiplicities later.

\begin{lemma}
Let $\lambda \vdash n$ be an arbitrary diagram, and let $\lambda'$
be its conjugate diagram. Consider two collections of diagrams: the
first one is the set of diagrams that are larger than (or equal to)
$\lambda$ in the dominance ordering; and the second one is the set
of diagrams conjugate to diagrams that are larger (in the same
sense) than $\lambda'$. The intersection of these two sets consists
of the single diagram $\lambda$.
\end{lemma}

The proof of this lemma consists in a direct check that
passing to the conjugate diagram inverts the dominance ordering.

\begin{lemma}
The irreducible representations
$\pi_{\mu}$ appearing in the induced representation
$\Pi_{\lambda}$ correspond to diagrams $\mu$ that are larger than
$\lambda$ in the dominance ordering: $\mu\unrhd\lambda$.
\end{lemma}

\begin{proof}
We will prove the lemma by induction, using the branching
rule for representations which we already have. The induction base
is, for example, the case
$n=3$, where the assertion is obvious. Assume that it holds for
$n-1$. Consider the representation
$\Pi_{\lambda}$ for some diagram
$\lambda\vdash n$. Its restriction to ${\frak
S}_{n-1}$ contain only those irreducible representations $\pi_{\rho}$
of ${\frak S}_{n-1}$ that appear as irreducible
components of the induced representations
$\Pi_{\gamma}$ corresponding to diagrams
$\gamma$ that precede $\lambda$, i.e., $\gamma\prec
\lambda$. Hence our assertion is reduced to the following combinatorial fact:
let $\lambda\vdash n$ and $\mu\vdash n$ be two diagrams; then the conditions

\noindent ($\alpha$) $\lambda\lhd\mu$

\noindent and

\noindent $(\beta)$ for every $\rho\vdash (n-1)$, $\rho \prec \mu$, there exists $\gamma\vdash
(n-1)$, $\gamma \prec \lambda$, such that $\rho \rhd \gamma$

\noindent are equivalent.

Obviously, $(\alpha)$ implies $(\beta)$. Now let $(\beta)$ hold and assume that the sum of the
lengths of a certain number $h$ of the first rows of $\lambda$ is greater than the same sum for
$\mu$ (i.e., $(\alpha)$ is not satisfied). Then, removing a cell from the row of $\mu$ with the
least possible number, we obtain a diagram $\rho\prec\mu$ for which condition $(\beta)$ is not
satisfied, because the sum of the same $h$ rows in any diagram $\gamma$ with  $\gamma\prec \lambda$
will be still greater than the sum of $h$ rows of $\rho$. Thus it follows from the above
considerations and the induction hypothesis that all irreducible components of the restriction of
the representation $\pi_{\mu}$ to the subgroup ${\frak S}_{n-1}$ appear in the restriction to the
same subgroup of the representation $\Pi_{\lambda}$, and hence the representation $\pi_{\mu}$ (with
 $\mu \rhd \lambda$) itself appears in the decomposition of $\Pi_{\lambda}$ into irreducible components.
 In the case $\mu=\lambda$ the condition above is also true so for the same reason as before
 $\pi_{\lambda}$ appears in $\Pi_{\lambda}$ as irreducible component. The same arguments show that
 if $\mu$ does not satisfy to the condition $\lambda \unlhd \mu$ then $\pi_{\mu}$ is not irreducible
 component of $\Pi_{\lambda}$.
\end{proof}

\begin{proof}[Proof of Theorem~1]
Let us proceed to the proof of Theorem~1.
Assume that the assertion is proved for the symmetric group of degree
$n-1$ and every induced representation
  $\Pi_{\gamma}$, $\gamma\vdash (n-1)$, of this group.
Consider an arbitrary diagram
$\lambda\vdash n$. Note that the set of diagrams $\gamma$ preceding
$\lambda$ is linearly ordered in the sense of the dominance ordering
and contains the minimal element, namely, the diagram that is obtained
from $\lambda$ by
removing the cell lying in the row with the least possible
number (the length of this row is the greatest possible); denote this diagram
by $\bar\lambda$ (since it depends only on  $\lambda$).
By the induction hypothesis, the representation $\pi_{\bar\lambda}$ appears in
the representation $\Pi_{\bar\lambda}$ with multiplicity one.
On the other hand, for any other
$\gamma$ distinct from $\bar\lambda$ and preceding
$\lambda$, it does not appear in the induced representation
$\Pi_{\gamma}$, because  $\bar\lambda$ cannot be larger in the
dominance ordering than any other diagram
$\gamma \prec \lambda$ except itself. But this means that
the multiplicity of $\pi_{\bar\lambda}$ in $Res_{{\frak
  S}_{n-1}}Ind_{H(\lambda)}^{{\frak S}_{n-1}}\textbf{1}$
is also equal to one. Thus the multiplicity of $\pi_{\lambda}$ in
the induced representation $\Pi_{\lambda}$ cannot be equal to zero
and cannot be greater than the multiplicity of any of the
restrictions of this representation to the subgroup ${\frak
S}_{n-1}$; therefore this multiplicity of the irreducible
representation $\pi_{\lambda}$ in the induced representation
$\Pi_{\lambda}$ is equal to one.

  Now, because the decomposition  of induced representation
    $Ind_{H(\alpha)}^{{\frak S}_n}\textbf{sgn}$
  differs from decomposition of induced representation
  $Ind_{H(\alpha)}^{{\frak S}_n}\textbf{1}$
  by exchange of irreducible components $\pi_{\mu}$
  onto representations $\pi_{\mu'}$ only, then
  irreducible representation $\pi_{\mu}$ belongs to the decomposition
  of the induced representation  $Ind_{V(\lambda)}^{{\frak S}_n}\textbf{sgn}$,
  (which is the same as $Ind_{H(\lambda')}^{{\frak S}_n} \textbf{sgn}$),
  iff for conjugate diagram to the diagram $\mu$ the following is true:
  $\mu' \unrhd \lambda'$. Thus we have $\mu\unlhd\lambda $, and therefore,
  accordingly to the lemma 3 the decomposition of the representation
  $Ind_{V(\lambda)}^{{\frak S}_n}\textbf{sgn}$, contains only
  representations  $\pi_{\mu}$ with $\mu\unlhd\lambda $, in
  particular representation $\pi_{\lambda}$; moreover - it has no
  multiplicity,
  and because of lemma 2 there are no other common components in the intersection
  of the decompositions of induced representations in the
  formulation of the theorem: $Ind_{H(\lambda)}^{{\frak S}_n}[\textbf{1}] \bigcap
  Ind_{V(\lambda)}^{{\frak S}_n}[\textbf{sgn}]=\{\pi\}$.
  \end{proof}

Denote the number of diagrams that majorize a given diagram
$\lambda$ in the dominance ordering by $h(\lambda)$,
and the number of diagrams that majorize
$\bar\lambda$ by ${\bar h}(\lambda)$. In what follows, we will need
these two numbers, as well as the operation $\lambda \rightsquigarrow \bar\lambda$
of removing a cell
from the uppermost possible row of $\lambda$.

\section{Corollaries: formulas for the multiplicities and
the recurrence property of Kostka numbers}

Above we have found out, by the inductive method, which irreducible
representations appear in the decompositions of induced representations
$\Pi_{\lambda}$. In order to obtain the complete decomposition, i.e.,
to find the multiplicities of irreducible components, we need
more detailed information on induced representations than that we have
used above. The formula for the multiplicities of the irreducible components
of induced representations is well known
(see \cite{JK, F, M}): the multiplicities are given by Kostka numbers (see below).
Usually, one obtains this formula by applying the theory of symmetric functions and
the theory of characters (the Frobenius formulas and so on), i.e.,
by, in a sense, nonelementary methods. Here we do not give a complete proof
of this formula, but

1) deduce from the previous formulas necessary recurrence conditions
on the multiplicities of irreducible components of induced representations of
two successive symmetric groups
 ${\frak S}_{n-1}$ and ${\frak S}_n$,

and

2) show that the Kostka numbers satisfy these conditions.

These formulas can be regarded as necessary conditions on
the multiplicities. Apparently, they exhaust all linear relations
between them. Sometimes they uniquely determine these
multiplicities by induction. In the general case, one should take into
account additional conditions satisfied by the multiplicities.

Recall that the Kostka number $K({\mu,\lambda})$ (see \cite[1.6]{Mc})
is the number of semistandard tableaux (strictly increasing in every column
and weakly increasing in every row) of shape
$\mu\vdash n$ and weight
$\lambda\vdash n$, which means that the number of symbols
$i$, $i=1,2,\dots, n$, in the weight is equal to the length of the $i$th
row of $\lambda$. It is clear that
$K({\mu,\lambda})> 0$ if and only if $\mu\unrhd\lambda$, and that
$K(\lambda,\lambda)=1$. For detailed information on
 $K({\mu,\lambda})$, see \cite{Mc}.

The complete decomposition of an induced representation, sometimes
called Young's rule, is as follows.

{\it Let $\lambda \vdash n$ be a diagram with $n$ cells. Then
 $$\Pi_{\lambda}\equiv Ind_{H(\lambda)}^{{\frak S}_n}\textbf{1}
 =\sum_{\mu:\mu \unrhd \lambda}K(\mu,\lambda)\pi_{\mu}=
 \pi_{\lambda}+\sum_{\mu:\mu \rhd
 \lambda}K(\mu,\lambda)\pi_{\mu},$$
where $\pi_{\mu}$ is the irreducible representation of
${{\frak S}_n}$ associated with a diagram $\mu$, and its multiplicity
is equal to the Kostka number $K(\mu,\lambda)$.}

Denote the multiplicity of the irreducible representation
$\pi_{\mu}$ in the induced representation
$\Pi_{\lambda}$ by $M(\mu,\lambda)$.

We will use the Frobenius--Young correspondence
$\lambda \Leftrightarrow \pi_{\lambda}$, assuming that it is already
established (by any method). Let
$\lambda \vdash n$, $\mu \vdash n$.

\begin{theorem}
Let $\lambda \vdash n$ and $\rho
\vdash (n-1)$ be arbitrary diagrams. Then
\begin{equation}
\sum_{\mu:\mu\succ\rho} M(\mu,
\lambda)=\sum_{\gamma:\gamma\succ\lambda}c(\lambda,\gamma)M(\rho,
\gamma),
 \end{equation}
where the sums range over diagrams
$\gamma\vdash(n-1)$ and $\mu\vdash n$.
\end{theorem}

\begin{proof}
Take the formula from Lemma~1 and, for a fixed
$\lambda$, decompose all induced representations of
${\frak S}_{n-1}$ in the right-hand side into irreducible components:
$$Res_{{\frak S}_{n-1}}Ind_{H(\lambda)}^{{\frak S}_n}
\textbf{1}=\sum_{\gamma:\gamma\prec\lambda}c(\lambda,\gamma)
\sum_{\rho:\rho\unrhd\gamma}M(\rho,\gamma)\pi_{\rho}.$$
Change the order of summation:
 $$=\sum_{\rho}\pi_{\rho}\sum_{\gamma:
\gamma\unrhd
\rho,\gamma\prec\lambda}c(\lambda,\gamma)M(\rho,\gamma).$$
The inner sum is the total multiplicity of the irreducible representation
$\pi_{\rho}$ in the restriction to ${\frak
S}_{n-1}$ of the induced representation $\Pi_{\lambda}$ of
${\frak S}_n$. On the other hand, it can be expressed in terms
of the desired multiplicities
$M(\mu,\lambda)$ by simply summing them over all irreducible representations
$\pi_{\mu}$ of ${\frak S}_n$ that appear in the decomposition of
$\Pi_{\lambda}$ and whose restriction to
${\frak S}_{n-1}$ contains $\pi_{\rho}$ as a subrepresentation:
$$\sum_{\gamma: \gamma\unlhd
\rho,\gamma\prec\lambda}c(\lambda,\gamma,)M(\rho,\gamma)
=\sum_{\mu:\mu \in \Lambda,\mu\succ\rho} M(\mu,\lambda). \eqno(*)$$

Since it follows from Theorem~1 that
$M(\mu,\lambda)\ne 0$ if and only if
$\mu\unrhd\lambda$, and the same holds for
$M(\rho,\gamma)$, and since these numbers are nonnegative, the formula
can be simplified:

\smallskip
for any $\lambda \vdash n$ and $\rho \vdash (n-1)$,
$$
\sum_{\mu:\mu\succ\rho} M(\mu,
\lambda)=\sum_{\gamma:\gamma\prec\lambda}c(\lambda,\gamma)M(\rho,
\gamma).
$$
\end{proof}

We leave the consideration of
the general case in the spirit of this approach
and a more detailed analysis of the problem for a suitable occasion.
Thus the multiplicities satisfy relations~(1). These relations are
a series of linear identities, and if we assume that the multiplicities
$M(\gamma,\rho)$ are found for all $\gamma,\rho$, i.e., for induced
representations of the group ${\frak S}_{n-1}$, then these identities
can be regarded as a system of linear equations on the multiplicities
$M(\mu,\lambda)$, i.e., for the group ${\frak
 S}_n$.

Now let us show that the Kostka numbers satisfy this system, i.e.,
relations~(1).

 \begin{theorem}
\begin{equation}
\sum_{\mu:\mu\succ\rho} K(\mu,
\lambda)=\sum_{\gamma:\gamma\prec\lambda}c(\lambda,\gamma)K(\rho,
\gamma).
 \end{equation}
\end{theorem}

\begin{proof}
We will give a ``bijective'' proof of this purely combinatorial fact.
In order to prove equality~(2), first let us understand that it
means counting the same number in two ways.
Namely, fix an arbitrary partition $\lambda$ of the set
$\{1, \ldots, n\}$ into blocks $(1^{r_1}, \ldots, k^{r_k},\ldots,
n^{r_n})$, $r_k\geq 0$, $k=1, \ldots, n$, $\sum_k kr_k=n$, $r_1\geq
r_2\geq \ldots \geq r_n\geq 0$, and a diagram $\rho$, $\rho \vdash
(n-1)$.

Then the number in question is the {\it number of semistandard
Young tableaux of shape $\rho$ with weights obtained from the weight
determined by $\lambda$ by removing one symbol.}
The left-hand side of the formula corresponds to counting this number
in the following order: first we find all diagrams
$\mu$ that majorize $\lambda$ in the dominance ordering and
are larger than $\rho$, and semistandard
tableaux of shapes $\mu$ and weight
$\lambda$, and then remove one cell from each of them so that to
obtain a semistandard tableau of shape
$\rho$. Conversely, the sum in the right-hand side corresponds to
counting the same number in another order; namely, first
we consider partitions (more exactly, weights)
$\gamma$ differing from
$\lambda$ by exactly one number (taking into account the multiplicities
$r_i$), and then enumerate all semistandard tableaux of shape
$\rho$ with the obtained weights $\gamma$. Obviously, in both cases
each semistandard tableau of shape
$\rho$ with one of the possible weights is obtained exactly once,
which proves the desired equality and, most importantly, determines
a {\it bijection, compatible with the natural ordering by inclusion and
the dominance ordering, between
the $\mu$-tableaux and $\rho$-tableaux}.
In fact, our formula asserts that two types of
operations are interchangeable: constructing a semistandard tableau of a given weight
on the one hand, and removing a cell from a diagram (in the left-hand
side) or removing an element from one of the blocks of a weight
(in the right-hand side).
\end{proof}

In the next section, we give an example, figures, and explanations to them.

\medskip
\noindent
 \textbf{Remark.}
 Formula (2) is a far generalization of the well-known relation
  $$n \cdot \dim (\rho) =\sum_{\mu:\mu \succ \rho} \dim (\mu).$$
   This relation is obtained if we take
   $\lambda$ to be a column, i.e.,
   $\lambda=1^n$. Then the corresponding weight is the collection of all numbers
   $\{1,2, \ldots, n\} $ without repetitions; thus
   $K(\lambda, \mu)=\dim (\mu)$ for all $\mu$, and the right-hand side is the sum
   of the dimensions of all diagrams larger than
   $\rho$. In the left-hand side, the weights obtained from
   $\lambda$ do not have repetitions either, and there
   are $n$ of them (since we successively remove each of the  elements of $\lambda$),
whence $c(\lambda,\gamma)=n$  and $K(\rho, \gamma)=\dim (\rho)$.
 \medskip

Relations~(2) should be regarded as a system of equations with respect to
$M(\mu,\lambda)$, and the right-hand sides should be assumed known. If we
subtract, term by term, system~(2) from system~(1), we will obtain
a homogeneous system with respect to the differences
$M(\mu,\lambda)-K(\mu,\lambda)\equiv Y_{\mu,\lambda}$:
\begin{equation}
\sum_{\mu:\mu\succ\rho} Y_{\mu,\lambda}=0.
\end{equation}

The parameter of the whole system~(3) is an arbitrary but fixed diagram
$\lambda$. If this system has only the zero solution, then the multiplicities
are equal to the Kostka numbers. Consider this system in more detail.
Above we have introduced the operation $\theta \rightarrow \bar \theta$ of removing the uppermost
possible cell from a diagram $\theta$,
denoted the number of diagrams that
majorize $\theta$ in the dominance ordering by
$h(\theta)$, and denoted $h(\bar\theta)$ by ${\bar h} (\theta)$. Note that
if $\mu\unrhd \lambda$ and $\mu\succ \rho$, then $\rho\unrhd \bar
\lambda$, and that the set of diagrams $\rho\vdash (n-1)$
that majorize one of the diagrams $\gamma \prec \lambda$ in the dominance ordering
coincides with the set of diagrams
$\rho$ that majorize $\bar
\lambda$. Note that if $\rho\unrhd\bar\lambda$, then adding
a cell to the row of $\rho$ with the least possible number gives
a diagram $\mu$ with $\mu\unrhd \lambda $, and, conversely, if $\mu\unrhd
\lambda $ and $\rho\prec\mu$, then $\rho\unrhd\bar\lambda$. Therefore
the equations of our system are indexed by the diagrams
$\rho$ satisfying the condition
$\rho\unrhd\bar\lambda$, and the unknowns are indexed by the diagrams
$\mu\unrhd\lambda$.\footnote{Though the conditions
$\rho\unrhd\bar\lambda$ and $\mu\unrhd\lambda$ in general do not imply that
$\mu\unrhd\lambda$.} In other words, the system is of order
${\bar h}(\lambda) \times h(\lambda)$.

\begin{statement} Assume that a diagram $\lambda$ is such that the mapping
 $\mu \rightarrow \bar\mu $ is a bijection between the set
 $\{\mu:\mu\unrhd\lambda\}$ and the set $\{\rho:\rho\unrhd\bar\lambda\}$.
 Then system~{\rm(3)} has only the zero solution, and
 thus the multiplicities are equal to the Kostka numbers.
\end{statement}

\begin{proof}
It follows from the assumptions that
 $h(\lambda) = {\bar h}(\lambda)$, so that the matrix of the system
 is a square matrix with entries 0 and 1.
It is not degenerate, because for
$\rho=\bar\mu\unlhd\mu$ and all $\rho'\prec \mu$ we have
 $\rho \unlhd \rho'$. By the same reason,
the matrix of the system is unipotent
with respect to the dominance ordering: if we identify the numbers of
unknowns and equations according to the bijection, then the entry
$a_{\rho,\rho'}$ of the matrix is zero unless
$\rho \unlhd \rho'$.
 \end{proof}

In the following example, this bijection condition is easy to check.

\begin{lemma}
If in a diagram $\lambda\vdash n$ the length of the first row is
at least as great as
the sum of the lengths of all the other rows (i.e.,
$\lambda_1 \geq n/2$), then the assumption of the previous proposition
is satisfied: the mapping $\mu
 \rightarrow \bar\mu $ is a bijection of the set $\{\mu:\mu\unrhd\lambda\}$
onto the image.
\end{lemma}

From the point of view of the inductive approach to the
representation theory of the symmetric groups, the  problem under
consideration can be formulated as follows: to what extent a representation of
${\frak S}_n$ is uniquely determined by its restriction to
${\frak S}_{n-1}$?  In general, there is no uniqueness unless we impose
 some additional conditions on the representation. In the case under
consideration, such a condition is that the representation is induced
from the identity representation of a Young subgroup,
and in order to conclude that the Kostka numbers are exactly the
required multiplicities, we need a uniqueness theorem that follows from additional (nonlinear)
identities for the multiplicities. This question is of interest
regardless of the fact that the multiplicities are known, because
its solution gives new relations between Kostka numbers. One may think that
the realization of induced representations in polylinear forms discussed below
will be useful in studying this problem.

Let us formulate another problem related to bijections between sets of
Young diagrams of two successive levels.

\smallskip\noindent
 \textbf{Problem.}
 The following question is close to the problems considered above,
 but it first arose in another context.

 Is it possible to define a natural correspondence (polymorphism)
 between the set ${\cal P}_{n-1}$ of diagrams with $n-1$ cells and the set
 ${\cal P}_n$ of diagrams with $n$ cells that would send the uniform
 distribution on the set ${\cal P}_{n-1}$ to the uniform
 distribution on the set ${\cal P}_n$? In other words, is it possible to define
 a nonnegative $p(n-1)\times p(n)$ matrix whose entries $c_{\gamma,\lambda}$,
 $\gamma \in {\cal P}_{n-1}$,
 $\lambda \in {\cal P}_n$, can be positive only if
 $\gamma \prec \lambda$ and whose rows (respectively, columns)
 sum to $p(n-1)^{-1}$ (respectively, $p(n)^{-1}$)?

\section{Examples}

We will give two examples illustrating the following key combinatorial fact
used in the first part of the theorem.
For every diagram $\lambda \vdash n$, the following two numbers coincide:

$\bullet$ the number of semistandard tableaux of weight
$\lambda$ with shapes $\mu$ satisfying the two conditions:

1) $\mu$ majorizes $\lambda$ in the dominance ordering, and

2) $\mu $ is larger (in the usual ordering) than a diagram
$\rho\vdash (n-1)$ (depending on $\mu$), where
$\rho$ majorizes one of the diagrams $\gamma\vdash (n-1) $ that are
less than $\lambda$;

and

$\bullet$ the number of semistandard tableaux of weight
$\lambda$ with one symbol removed with shapes
$\rho$ satisfying the two conditions:

1) $\rho$ majorizes, in the dominance ordering,  one of the diagrams
$\gamma$ that are less (in the usual ordering) than $\lambda$, and

2) $\rho$ is less (in the usual ordering) than a diagram
$\mu$ that majorizes $\lambda$ in the dominance ordering.

We will describe the canonical bijection between the semistandard
$\mu$- and $\rho$-tableaux. Schematically, the situation is represented by
the following figure.

Diagrams $\lambda$ and $\rho$ are fixed; diagrams
$\gamma$ and $\mu$ and tableaux of shape
$\rho$ vary. We will show how to establish a bijection between the set of
transitions from $\lambda$ to  $\rho$ through $\mu$
and the set of transitions from $\lambda$ to  $\rho$ through
$\gamma$, i.e., between the sets of two types of paths in the following scheme:

\def\circrho{\hbox{$\bigcirc$ \kern-1.15em $\rho$}}
\def\circlambda{\hbox{$\bigcirc$ \kern-1.15em $\lambda$}}

{\LARGE$$\begin{matrix}
\boxit\gamma   &\unlhd & \circrho\\
\curlywedge &    \qquad   & \curlywedge \\
\circlambda & \unlhd &  \boxit\mu\\
\end{matrix}$$}

Here $\unlhd$ stands for the dominance ordering,  and
$\prec$ stands for the usual ordering of diagrams by inclusion. Transitions (paths)
of the first type are determined by a semistandard tableau of shape
$\mu$, and transitions of the second type are determined  by a semistandard tableau of shape
$\rho$, since the diagram
$\gamma$ is uniquely determined by the
$\rho$-tableau. Hence a bijection between the paths is determined by a
correspondence between the semistandard
$\mu$-tableaux and semistandard $\rho$-tableaux.

In the first example, the weight has no multiplicities (i.e.,
$\lambda$ has no rows of equal lengths), and in the second example,
there are multiplicities; however, the algorithm for establishing the
required bijection does not essentially depend on this fact.

\bigskip
\noindent{\bf Example 1.} Let
$n=6$ and $\lambda=(3,2,1) =\,\lower2\cellsize\vbox{\footnotesize
\cput(1,1){1}
\cput(1,2){1}
\cput(1,3){1}
\cput(2,1){2}
\cput(2,2){2}
\cput(3,1){3}
 \cells{
 _ _ _
|_|_|_|
|_|_|
|_|}}$\,.

\smallskip
Below we successively list the
$\mu$-tableaux, $\gamma$-weights, and $\rho$-tableaux.

\bigskip

$\mu:A=\,\lower1\cellsize\vbox{\footnotesize
\cput(1,1){1}
\cput(1,2){1}
\cput(1,3){1}
\cput(1,4){2}
\cput(1,5){2}
\cput(2,1){3}
 \cells{
 _ _ _ _ _
|_|_|_|_|_|
|_|}}$ ,
$B=\,\lower1\cellsize\vbox{\footnotesize
\cput(1,1){1}
\cput(1,2){1}
\cput(1,3){1}
\cput(1,4){2}
\cput(1,5){3}
\cput(2,1){2}
 \cells{
 _ _ _ _ _
|_|_|_|_|_|
|_|}}$ ,
$C=\,\lower1\cellsize\vbox{\footnotesize
\cput(1,1){1}
\cput(1,2){1}
\cput(1,3){1}
\cput(1,4){2}
\cput(2,1){2}
\cput(2,2){3}
 \cells{
 _ _ _ _
|_|_|_|_|
|_|_|}}$ ,
$D=\,\lower1\cellsize\vbox{\footnotesize
\cput(1,1){1}
\cput(1,2){1}
\cput(1,3){1}
\cput(1,4){3}
\cput(2,1){2}
\cput(2,2){2}
 \cells{
 _ _ _ _
|_|_|_|_|
|_|_|}}$ ,
$E=\,\lower1\cellsize\vbox{\footnotesize
\cput(1,1){1}
\cput(1,2){1}
\cput(1,3){1}
\cput(1,4){2}
\cput(2,1){2}
\cput(2,2){3}
\cells{
 _ _ _ _
|_|_|_|_|
|_|_|}}$;
\bigskip

$\gamma:X=\,\lower2\cellsize\vbox{\footnotesize
\cput(1,1){1}
\cput(1,2){1}
\cput(2,1){2}
\cput(2,2){2}
\cput(3,1){3}
 \cells{
 _ _
|_|_|
|_|_|
|_|}}$ ,
$Y=\,\lower2\cellsize\vbox{\footnotesize
\cput(1,1){1}
\cput(1,2){1}
\cput(1,3){1}
\cput(2,1){2}
\cput(3,1){3}
 \cells{
 _ _ _
|_|_|_|
|_|
|_|}}$ ,
$Z=\,\lower1\cellsize\vbox{\footnotesize
\cput(1,1){1}
\cput(1,2){1}
\cput(1,3){1}
\cput(2,1){2}
\cput(2,2){2}
 \cells{
 _ _ _
|_|_|_|
|_|_|}}$ ;
\bigskip

$\rho:L=\,\lower1\cellsize\vbox{\footnotesize
\cput(1,1){1}
\cput(1,2){1}
\cput(1,3){2}
\cput(1,4){2}
\cput(2,1){3}
 \cells{
 _ _ _ _
|_|_|_|_|
|_|}}$ ,
$M=\,\lower1\cellsize\vbox{\footnotesize
\cput(1,1){1}
\cput(1,2){1}
\cput(1,3){2}
\cput(1,4){3}
\cput(2,1){2}
 \cells{
 _ _ _ _
|_|_|_|_|
|_|}}$ ,
$N=\,\lower1\cellsize\vbox{\footnotesize
\cput(1,1){1}
\cput(1,2){1}
\cput(1,3){1}
\cput(1,4){2}
\cput(2,1){3}
 \cells{
 _ _ _ _
|_|_|_|_|
|_|}}$ ,
$P=\,\lower1\cellsize\vbox{\footnotesize
\cput(1,1){1}
\cput(1,2){1}
\cput(1,3){1}
\cput(1,4){3}
\cput(2,1){2}
 \cells{
 _ _ _ _
|_|_|_|_|
|_|}}$ ,
$Q=\,\lower1\cellsize\vbox{\footnotesize
\cput(1,1){1}
\cput(1,2){1}
\cput(1,3){1}
\cput(1,4){2}
\cput(2,1){2}
 \cells{
 _ _ _ _
|_|_|_|_|
|_|}}$ .
\bigskip

Transition from $\gamma$ to a $\rho$-tableau:
$$X\longrightarrow L,M; \quad Y\longrightarrow N,P;\quad Z\longrightarrow Q.$$

Bijection $\mu \longleftrightarrow \rho$:

$$A\longleftrightarrow L;\quad B\longleftrightarrow M;\quad C\longleftrightarrow N;\quad D\longleftrightarrow P;\quad  E\longleftrightarrow Q.$$

\bigskip
\noindent
{\bf Example 2.} In this example, the weight has multiplicities, that is,
the coefficient $c$ (see the proof) is not always equal to $1$, so that the bijection
is established in a more general way.

Let $n=5$ and $\lambda=(2^2,1) =\,\lower2\cellsize\vbox{\footnotesize
\cput(1,1){1}
\cput(1,2){1}
\cput(2,2){2}
\cput(2,1){2}
\cput(3,1){3}
 \cells{
 _ _
|_|_|
|_|_|
|_|}}$ . Then

\bigskip
\bigskip

$\mu:A=
\,\lower1\cellsize\vbox{\footnotesize
\cput(1,1){1}
\cput(1,2){1}
\cput(1,3){2}
\cput(1,4){2}
\cput(2,1){3}
 \cells{
 _ _ _ _
|_|_|_|_|
|_|}}$ ,
$B=
\,\lower1\cellsize\vbox{\footnotesize
\cput(1,1){1}
\cput(1,2){1}
\cput(1,3){2}
\cput(1,4){3}
\cput(2,1){2}
 \cells{
 _ _ _ _
|_|_|_|_|
|_|}}$ ,
$C=
\,\lower1\cellsize\vbox{\footnotesize
\cput(1,1){1}
\cput(1,2){1}
\cput(1,3){2}
\cput(2,1){2}
\cput(2,2){3}
 \cells{
 _ _ _
|_|_|_|
|_|_|}}$ ,
$D=
\,\lower1\cellsize\vbox{\footnotesize
\cput(1,1){1}
\cput(1,2){1}
\cput(1,3){3}
\cput(2,1){2}
\cput(2,2){2}
 \cells{
 _ _ _
|_|_|_|
|_|_|}}$ ,
$E=
\,\lower2\cellsize\vbox{\footnotesize
\cput(1,1){1}
\cput(1,2){1}
\cput(1,3){2}
\cput(2,1){2}
\cput(3,1){3}
 \cells{
 _ _ _
|_|_|_|
|_|
|_|}}$ ;

\bigskip
$\gamma: X=
\,\lower2\cellsize\vbox{\footnotesize
\cput(1,1){1}
\cput(1,2){2}
\cput(2,1){2}
\cput(3,1){3}
 \cells{
 _ _
|_|_|
|_|
|_|}}$ ,
$Y=
\,\lower2\cellsize\vbox{\footnotesize
\cput(1,1){1}
\cput(1,2){1}
\cput(2,1){2}
\cput(3,1){3}
 \cells{
 _ _
|_|_|
|_|
|_|}}$ ,
$Z=
\,\lower1\cellsize\vbox{\footnotesize
\cput(1,1){1}
\cput(1,2){1}
\cput(2,1){2}
\cput(2,2){2}
 \cells{
 _ _
|_|_|
|_|_|}}$ ;

\bigskip

$\rho:L=
\,\lower1\cellsize\vbox{\footnotesize
\cput(1,1){1}
\cput(1,2){2}
\cput(1,3){2}
\cput(2,1){3}
 \cells{
 _ _ _
|_|_|_|
|_|}}$ ,
$M=
\,\lower1\cellsize\vbox{\footnotesize
\cput(1,1){1}
\cput(1,2){2}
\cput(1,3){3}
\cput(2,1){2}
 \cells{
 _ _ _
|_|_|_|
|_|}}$ ,
$N=
\,\lower1\cellsize\vbox{\footnotesize
\cput(1,1){1}
\cput(1,2){1}
\cput(1,3){2}
\cput(2,1){3}
 \cells{
 _ _ _
|_|_|_|
|_|
}}$ ,
$P=
\,\lower1\cellsize\vbox{\footnotesize
\cput(1,1){1}
\cput(1,2){1}
\cput(1,3){3}
\cput(2,1){2}
 \cells{
 _ _ _
|_|_|_|
|_|}}$ ,
$Q=
\,\lower1\cellsize\vbox{\footnotesize
\cput(1,1){1}
\cput(1,2){1}
\cput(1,3){2}
\cput(2,1){2}
 \cells{
 _ _ _
|_|_|_|
|_|}}$ .

\bigskip
Transition from $\gamma$ to $\rho$:
$$X\longrightarrow L,M; \quad Y\longrightarrow N,P;\quad Z\longrightarrow Q.$$

Bijection $\mu \longleftrightarrow\rho$ :
$$A\longleftrightarrow L; \quad B\longleftrightarrow M; \quad  C\longleftrightarrow N; \quad D\longleftrightarrow P; \quad E\longleftrightarrow Q.$$

\bigskip
We will explain only the second example, which is slightly more complicated.

We are given two diagrams $\lambda \vdash 5$ and $\rho \vdash 4$;
the diagram $\lambda$ is regarded both as a diagram and a weight
$(1^2,2^2,3^1)$; the semistandard tableaux of shape
$\rho$ have the weights $\gamma$ that  are obtained from the weight
$\lambda$ by removing one of the symbols
(1, or 2, or 3). In the last row, we have listed all such fillings of
the diagram $\rho$. They can be obtained in two ways:

(1) The first way. We should consider all diagrams
$\mu$ that majorize $\lambda$ in the dominance ordering and
are larger in the usual ordering than
a diagram $\rho$ that has one cell less than
$\mu$; there are three such diagrams: (4,1), (3,2), (3,1,1).
There are five semistandard tableaux of weight $\lambda$.
These diagrams and semistandard tableaux are listed in the first row:
$A$, $B$, $C$, $D$, $E$; two diagrams, (4,1) and (3,2), allow
two semistandard tableaux each,
$A$, $B$ for the first one, and $C$, $D$ for the second one,
and the last diagram allows one tableau $E$. In the first two cases,
$K(\mu,\lambda)=2$, and in the third one, $K(\mu,\lambda)=1$.
The right-hand side of our formula equals $5$. Then we should remove the symbol
$1$ from the tableaux $A$ and $B$, the symbol
$2$ from the tableaux $C$ and $D$, and the symbol $3$ from the tableau $E$,
so that to obtain tableaux of shape $\rho$; again, there are five such tableaux.

The correspondence between the semistandard tableaux of shapes
$\mu$ and semistandard tableaux of shapes $\rho$ is given above.

(2) The second way. We should remove one of the cells of
$\lambda$, obtaining semistandard tableaux of different shapes
$\gamma$; there are three of them: $X$, $Y$, $Z$. But since
$\lambda$ has a multiple row
(of length $2$), it follows that for the diagram $\gamma=(2,1^2)$
(the tableaux $X$, $Y$), the coefficient $c(\lambda, \gamma)$ equals
two, and for $\gamma=2^2$ (the tableau $Z$),
this coefficient equals one. Then, regarding
$\gamma$ as a weight, we should fill the diagram
$\rho$, thus obtaining a semistandard tableau of shape
$\rho$. The correspondence between weights and tableaux is as follows:
the diagrams $X$, $Y$ generate two tableaux each, and
$Z$ generates one tableau $Q$.

\section{Realization of induced representations in polylinear forms}

\subsection{Definition of the space of forms}

Now let us consider induced representations from another point of view; we will describe their
realization in polylinear forms. Various particular cases of such realizations were considered
earlier. We keep the previous notation: $\lambda\vdash n$ is an arbitrary diagram, which determines
the type of a partition of the set $\{1, \dots ,n\}$ and the representation $
\Pi_{\lambda}=Ind_{H(\lambda)}^{{\frak S}_n} \textbf{1}$ induced from the identity representation
of the Young subgroup $H(\lambda)$; and $\pi_{\lambda}$ is the irreducible representation
corresponding to the diagram $\lambda$. In what follows, we will consider $q$-ary polylinear forms
with complex coefficients  in $n$ commuting variables, for some $q$ depending on the diagram. Of
course, the language of polylinear forms can be replaced by the language of tensors, but forms are
more convenient for our purposes.

Denote by $\lambda_1,\lambda_2, \dots,\lambda_k$, $\sum_{i=1}^k \lambda_i =n$,
the lengths of the rows of the diagram
$\lambda$, and by $n'=n-\lambda_1$ the number of cells in the rows
starting from the second one. Consider the following monomial in
$n'$ variables $y_{i,j}$ indexed by the cells
$(i,j)$ of $\lambda$ lying in the rows starting from the second one
 $(i>1)$:
  $$X_{\lambda}(\{y_{i,j}\}) \equiv \prod_{i,j}{[y_{i,j}]}^{i-1},$$
 The degree (``arity'') of this monomial is equal to $\sum_{j=1}^m C_{r_j}^2$, where $m$ is the number
of columns of $\lambda$ and $r_j$ is the length of the $j$th column.

Consider the complex linear span
$L_{\lambda}$ of all monomials of the form
 $X_{\lambda}$ in $n$ variables $x_1, \dots, x_n$ (i.e., we take as
$\{y_{i,j}\}$ all possible sets of $n'$ variables chosen from $x_1,\dots, x_n$).
This is a space of complex polylinear forms
in $n$ variables, which is obviously invariant under the
action of the group
${\frak S}_n$ by substitutions of variables in the space of all complex
polylinear forms.

 \begin{statement}
The constructed representation of ${\frak S}_n$ in the space $L_{\lambda}$
is equivalent to the induced representation
 $\Pi_{\lambda}$.
 \end{statement}

\begin{proof}
Choose an ordered collection of $n'$ variables $x_{s_1},\dots,
x_{s_{n'}}, n'=n-\lambda_1$ and consider the monomial
$X_{\lambda}(x_{s_1},\dots, x_{s_{n'}})$ in this variables. Identify
this monomial with a partition of type
 $\lambda$ of the set of variables $x_1,\dots, x_n$ as follows:
put the variables occurring in the monomial into the cells of $\lambda$ starting from the second
row so that the degree of the variable put into the cell $(i,j)$ in the monomial be equal to $i-1,
i>1$, $j=1,\dots, k$, and the cells of the first row (with $i=1$) fill by the variables
 $x_s$ with $s\ne s_l$, $l=1, \dots, s_{n'}$ (i.e., occur with zero exponent
in the monomial) in an arbitrary way.

By the definition of the monomials $X_\lambda$, it is always possible to do this, and we obtain a
partition of the set of $n$ variables according to the rows of the diagram $\lambda$. This gives a
bijection between the set of all partitions of type ${\lambda}$ and the set of monomials of type
$X_{\lambda}$ (recall that the number of the variable in the blocks are increased, the blocks of
partitions are listed in decreasing order of their lengths, and partitions with distinct numeration
of blocks of the same length are assumed distinct). It is also obvious that the constructed mapping
that sends monomials to partitions is bijective. Hence the number of monomials, as well as the
number of partitions, is equal to the multinomial coefficient
$\frac{n}{\prod_{i=1}^k(\lambda_i!)}$.

Thus the set of monomials is identified with the homogeneous space
with respect to the Young subgroup $H(\lambda)$, and the space  $L_{\lambda}$
of forms of the given
type is identified with the space of functions on the homogeneous space,
so that it is the space of the induced representation
$\Pi_{\lambda}$.
\end{proof}

The problem of finding an {\it explicit} decomposition of the
representation of the symmetric group induced from a Young subgroup
in the general case is very difficult; moreover,
an invariant problem is to decompose it into primary rather than irreducible
components, and, as far as the author knows, it has not yet an explicit
solution except in a number of particular cases.
The suggested realization helps to analyze the decompositions, because
the spaces of forms possess additional structures.

\subsection{Examples}

\noindent{\bf1.} For the row diagram $\lambda$ (i.e., the trivial partition),
the Young subgroup is the whole symmetric group; the unique, up to
a factor, monomial is a scalar, and the space
$L_{\lambda}$ is one-dimensional; we obtain the one-dimensional
identity representation.

\medskip
\noindent{\bf2.} On the other hand, the column diagram corresponds to the partition
of  $\{1,2,\dots, n\}$ into $n$ separate points; the Young subgroup is
the identity subgroup, and hence the induced representation is the
regular representation. The monomials
$X_{\lambda}$ in this case have the form
 $y_1^{n-1}y_2^{n-2}\dots y_{n-1}$ ($n-1$ variables);
there are $n!$ of them, and
 their degree equals
$C_n^2$. A monomial
 $x_{i_1}^{n-1}x_{i_2}^{n-2}\dots x_{i_{n-1}}^{n-1}$  is identified
with an ordered partition into separate points, i.e., with a
substitution
 $s\vdash i_s$, $s=2,\dots, n$ (the image of $x_1$
is the only variable that does not occur in the monomial).
The decomposition problem in this case is the problem of decomposing
the regular representation, which is an evidence that the general problem
is complicated.

\medskip
\noindent{\bf 3.} Let us consider in more detail the case of two-block partitions,
 i.e., two-row Young diagrams, in terms of forms. The Young subgroup
 in this case is the product of two symmetric groups.

Fix $n$ and $k$, and thus the diagram $\lambda=(n-k,k)$, $k \leq
[n/2]$. Let us show that  $\Pi_{\lambda}$ can be realized as a substitutional
representation in the space
$F_k$ of all square-free $k$-ary forms in $n$ variables, where the group
${\frak S}_n$ acts by substitutions of variables.

Consider the linear space $L_{t,k}\subset F_k$
generated by the $k$-forms $X_I$, where $I=\{i_1,\dots, i_{2t}\}$
runs over all ordered collections of
$2t$ indices from  $1, \dots, n$ and
$$X_I=(x_{i_1}-x_{i_2})(x_{i_3}-x_{i_4})\dots (x_{i_{2t-1}}-x_{i_{2t}})
\sigma_{k-t}(x_{j_1}x_{j_2} \dots x_{j_{n-2t}});\ j_l \in \{1,\dots,
n\}\setminus I,\; l= 1,\dots, n-2t;
$$
here $\sigma_p(\dots)$ is the basic
symmetric function of degree
$p$, i.e., the sum of the products of $p$-collections of distinct
variables over all such collections. By definition, the spaces
$L_{t,k}$ are invariant under the action of the group.

The following assertion can be checked directly:

\begin{statement}
{\rm1)} The representation $\pi_{n-l,l}$ is equivalent to the representation
in the subspace $L_{t,l}$.

{\rm 2)} We have the decomposition $$F_k=\oplus_{l=0}^k L_{t,l},$$
which coincides with the decomposition into irreducible components of the induced
representation $\Pi_{\lambda}$:
 $$\Pi_{\lambda}=\sum_{l=0}^k \pi_{n-l,l}.$$
The representation $\Pi_{\lambda}$ is multiplicity-free.

{\rm 3)} For even $n$ and $k=n/2$, the highest component, corresponding to
$l=k=n/2$ (i.e.,
the subspace $L_{n/2,n}$) consists of all $k$-forms depending on
the differences of variables.
\end{statement}

Note that the factor $\sigma_p(\cdot)$ is not needed for the realization of one irreducible
representation $\pi_{n-l,l}$ taken individually, and in fact the last claim of the proposition says
that every representation $\pi_{n-k,k}$ can be realized in the space of square-free $k$-forms
depending on the differences of variables (see below and \cite{TV}). But the factor
$\sigma_p(\cdot)$ is needed to show how the lower irreducible components are embedded into the
induced representation, i.e., the space of forms. Precisely this embedding of irreducible (or
primary) representations into the induced representation is the goal of describing the
decompositions for an arbitrary Young subgroup. The main difficulty is that the decompositions of
induced representations corresponding to diagrams with more than two rows involve multiplicities.

As concerns two-row diagrams themselves, the decomposition of the induced representation into
irreducible components can be described (and this description is more convenient) in terms of
decompositions of spaces of symmetric tensors (see \cite{VT, Nik}). The absence  of multiplicities,
for all $k$ the representations $\pi_{n-k,k}$ are embedded into the representation $\pi_{n,[n/2]}$,
i.e., into the space of tensors of valence $[n/2]$ and dimension $n$. However, the language of
tensors is not as convenient in more complicated cases of general diagrams; here it makes sense to
use the less common language of forms, which should be regarded as an extension of the tensor
techniques.

The above decomposition can easily be reformulated in terms of partitions,
i.e., in terms of the original description of the induced representation,
which acts,
by definition, in the space of functions on the space of partitions, which is
the homogeneous space with respect to the Young subgroup; for this
it suffices to observe that in the case of partitions into two sets,
we can identify a partition with one (the smaller) of them, and then
identify this set, for example,
$(i_1,i_2,\dots,  i_t)$, with the form (monomial)
$x_{i_1}\cdot x_{i_2} \ldots  \cdot x_{i_t}$.

\medskip
\noindent{\bf4.}
The next, in order of difficulty, example, after the two-row one,
which already involves multiplicities and has some features of the general situation,
is the partition $\lambda=(2,1^2)$ of the set of four elements. The corresponding
Young group is a group of order two. The decomposition of the induced
representation contains four irreducible representations
$\pi_{\mu}$, where $\mu$ runs over the diagrams
$1^4$, $(2,2)$, $(2,1^2)$, $(3,1)$ of dimensions
$1,2,3,3$, respectively, the last one with multiplicity $2$. For brevity, denote these
representations by $\pi_i$, $i=1,2,3,4$. In this case, the monomial
has the form
$$X_{\lambda}(y,z)=y^2\cdot z$$
and depends on two variables (according to the number of cells in the
2nd and 3rd rows of the diagram  $\lambda=(2,1^2)$); elements of the 12-dimensional space
$L_{\lambda}$ are forms depending on four variables
$x_i$, $i=1,2,3,4$.

Let us describe the decomposition of the space
$L_{\lambda}$ into the subspaces of irreducible representations,
indicating a convenient basis in each of them:

$\pi_1$ is a one-dimensional subspace, a basis consists of the form
$$D(x_1,\dots x_4)=\sum_{i\ne
  j}x_i^2x_j;$$

  $\pi_2$ is a two-dimensional subspace, a basis is given by
  \begin{eqnarray*}C_1(x_1,\dots x_4)&=&(x_1-x_2)(x_3-x_4)(x_1+x_2+x_3+x_4),\\
  C_2(x_1,\dots x_4)&=&(x_1-x_3)(x_2-x_4)(x_1+x_2+x_3+x_4)
  \end{eqnarray*}
(the form $C_3=(x_1-x_4)(x_2-x_3)(x_1+x_2+x_3+x_4)$ is a linear combination of
$C_1,C_2$);

  $\pi_3$ is a three-dimensional subspace (the ``main representation''
  with diagram  $\lambda$); a basis consists of the Specht polynomials
  (see below)
  \begin{eqnarray*}
  SP_1&=&(x_1-x_2)(x_1-x_3)(x_2-x_3),\\
  SP_2&=&(x_2-x_3)(x_2-x_4)(x_3-x_4),\\
  SP_3&=&(x_1-x_3)(x_1-x_4)(x_3-x_4);
  \end{eqnarray*}
the other Specht polynomials are linear combinations of
 $SP_1,SP_2,SP_3$;

 $\pi_4$ is a three-dimensional (natural) representation having
 multiplicity $2$ in the induced representation under consideration.
 Remarkably, we can canonically (in terms of forms) pick out two invariant
 subspaces  $A$ and
 $B$, each of them corresponding to the irreducible representation
  $\pi_4$.

 The first representation, $\pi_4^A$, has a basis
 $$A_k=\sum_{i\ne j}\epsilon_{i,j}^k x_ix_j, \qquad k=1,2,3,$$
 where
 $$
 \epsilon_{i,j}^k=\begin{cases}
 +1 &\mbox{if } i=k \mbox{ or } j=k,\\
 -1 & \mbox{otherwise},
 \end{cases}\qquad k=1,2,3.$$

The second one, $\pi_4^B$, has a basis
\begin{eqnarray*}
B_1&=&x_1^2(x_2+x_3+x_4)-x_1(x_2^2+x_3^2+x_4^2),\\
B_2&=&x_2^2(x_1+x_3+x_4)-x_2(x_1^2+x_3^2+x_4^2),\\
B_3&=&x_3^2(x_1+x_2+x_4)-x_3(x_1^2+x_2^2+x_4^2).
\end{eqnarray*}

The invariant way to distinguish between these subspaces
of representations equivalent to $\pi_4$ is
based on the fact that the representation
$\pi_4^A$, as well as the representations $\pi_1,\pi_2$, acts
in the six-dimensional subspace of {\it even forms}, while
the representation $\pi_4^B$, as well as $\pi_3$, acts
in the six-dimensional subspace of {\it odd forms}. Note that
the difficulty of the problem of decomposing representations into
irreducible components lies precisely in the fact that in the case of
multiplicities there is no natural decomposition into irreducible
components, but only into primary ones; but here (and, apparently, in the
general case of realizations in spaces of forms), the structure of the space of forms
allows one to obtain even an invariant decomposition. The point is that
there is another group acting on the space of forms, the group
of ``substitutions of degrees'';
here this group is ${\Bbb Z}_2$. Since in this example monomials involve
only two variables, it essentially gives a detailed
decomposition of the representation in square-free bilinear forms, i.e.,
in the space of 2-tensors, or matrices with zero diagonal; even (odd)
forms correspond to symmetric (skew-symmetric) matrices; but even the form
of the basis shows the advantage of the language of forms.

  \subsection{Specht modules and Specht polynomials}

Consider the so-called Specht polynomial $SP_t$ corresponding to an arbitrary  tableau
$t$ of shape $\lambda$ filled by variables
$x_{i,j}$ from the ground field, where
the numbers of variables correspond to the numbers of cells in the tableau $t$:
            $$SP_T=\prod_j \prod_{s<k} (x_{s,j}-x_{k,j}); $$
here the indices $s,k$ in each factor
range over the numbers of variables in the $j$th column. Thus
        $SP_T$ is the product of the Vandermonde determinants over all columns.
        The degree of $SP_T$ is equal to the sum of the numbers of pairs
        of variables in the columns; the polynomial does not contain the variable
        occupying a cell of the first row if the corresponding column consists
        only of this cell. The linear span of the functions obtained from
        $SP_T$ by permutations of all variables
        is the vector space of a representation
        corresponding to the diagram
        $\lambda$. It is called the Specht module $SP_{\lambda}$.

This definition implies

 \begin{statement}
 Every Specht polynomial (regarded as a form) is an element of the space of forms
 $L_{\lambda}$ defined above: $SP_{\lambda}\subset L_{\lambda}$.
 \end{statement}

 Indeed, it is not difficult to check that all monomials appearing in
 the product of the Vandermonde determinants have the form
 of monomials $X_{\lambda}$.
Evidently Specht polynomials depend only on the differences of variables, more exactly, they are
invariant under a shift, i.e., under simultaneously adding a constant to all variables. It is not
difficult to prove also the converse.

 Let us give basic facts on the Specht modules.

\begin{theorem}
{\rm 1)} The Specht module $SP_{\lambda}$ is the submodule of the space
of forms $L_{\lambda}$ defined above that consists of all forms depending
on the differences of variables.

 {\rm 2)} The Specht module $SP_{\lambda}$, as a representation of ${\frak
S}_n$, coincides with the irreducible representation
$\pi_{ \lambda}$ corresponding to the diagram $\lambda$.

 {\rm 3)} Denote by $SP_t$ the Specht polynomial corresponding to the arrangement of
 variables according to a standard tableau
  $t$. The polynomials  $SP_t$, where $t$ ranges over all standard
 tableaux of shape $\lambda$,
  form a basis of the Specht module $SP_{\lambda}$.
\end{theorem}

Since the irreducible representation $\pi_{ \lambda}$
occurs with multiplicity one in the induced representation
$\Pi_{\lambda} =Ind_{H(\lambda)}^{{\frak S}_n}\textbf{1}$
(the Frobenius--Young correspondence),
it follows that this unique subrepresentation is precisely the Specht module.

The proof of claim~1) follows from definitions. The most important part of this theorem, claims~2)
and~3), is known in this formulation; usually, one proves it using the theory of characters and
symmetric functions (see \cite{Mc,M,F}). We mean another proof, based on the same ideas that were
used in this paper; we will consider it in detail elsewhere.


\begin{thebibliography}{15}

\bibitem{Curt} C.~W.~Curtis and I.~Reiner, {\it
Representation Theory of Finite Groups and Associative Algebras},
Interscience Publ., New York--London, 1962. Russian translation:
Nauka, Moscow, 1969.
 \bibitem{Fu}W.~Fulton, {\it Young Tableaux. With Applications to
 Representation Theory and Geometry}. Cambridge University Press,
 Cambridge, 1997.
 \bibitem{F} W.~Fulton and J.~Harris, {\it Representation Theory.
 A First Course}. Springer-Verlag, New York, 1991.
\bibitem{Ham} M.~Hamermesh, {\it Group Theory and Its Application to Physical
Problems}. Addison-Wesley, Reading, Mass.--London, 1962.
Russian translation: Mir, Moscow, 1966.
\bibitem{James} G.~James, {\it The Representation Theory
of the Symmetric Group}. Springer, Berlin, 1978.
\bibitem{JK} G.~James and A.~Kerber,  {\it The Representation Theory
 of the Symmetric Group}. Addison-Wesley, Reading, Mass., 1981.
\bibitem{Mc} I.~Macdonald, {\it Symmetric Functions and Hall Polynomials},
2nd edition.
Clarendon Press, Oxford, 1995.
Russian translation of the first edition: Mir, Moscow, 1985.
\bibitem{Ma}G.~Mackey, {\it The Theory of Unitary Group Representations}.
 The University of Chicago Press, Chicago--London, 1976.
 \bibitem{M}C.~Musili, {\it Representations of Finite Groups}.
 Hindustan Book Agency, Delhi, 1993.
\bibitem{Nik}P.~P.~Nikitin,  A realization of the
irreducible representations of ${\frak S}_n$
corresponding to 2-row diagrams in
square-free symmetric multilinear forms, {\it Zap. Nauchn. Semin. POMI}
{\bf 301} (2003), 212--219. English translation:
{\it J. Math. Sci. (New York)}  {\bf129}, No. 2 (2005), 3796--3799.
\bibitem{OV} A.~Okounkov and A.~Vershik, New approach to representation
 theory of symmetric groups, {\it
 Selecta Math.} {\bf 2}, No.~4 (1996), 1--15.
\bibitem{Se}J.-P.~Serre, {\it Linear Representations of Finite
Groups}, Springer-Verlag, New York--Heidelberg, 1977. Russian translation:
Mir, Moscow, 1970.
\bibitem{TV}
 N.~V.~Tsilevich and A.~M.~Vershik, On different models
 of representations of the infinite symmetric group.
 To appear in {\it Adv. Appl. Math.}
 \bibitem{V} A.~M.~Vershik, Local algebras and a new version of
 Young's orthogonal form, In: ``Topics in Algebra'', Banach Cent. Publ.,
 V.~26. Part~2. PWN-Polish Sci.
 Publ. Warszawa (1990), pp.~467--473.
\bibitem{VO}A.~M.~Vershik and A.~Yu.~Okounkov,
A new approach to the representation theory of symmetric groups. II.
{\it Zapiski Nauchn. Semin. POMI} {\bf 307} (2004), 57--98.
English translation: {\it J. Math. Sci. (New York)} {\bf 131}, No.~2 (2005),
5471--5494.
\bibitem{VT} A.~M.~Vershik and N.~V.~Tsilevich,
Markov measures on Young tableaux
and induced representations of the infinite symmetric group,
{\it Probab. Theory Appl.} {\bf51} (2006), 47--63.
\bibitem{Weyl} H.~Weyl, {\it The Theory of Groups and Quantum Mechanics},
Dover Publ., New York, 1949. Russian translation:
Nauka, Moscow, 1986.
\end{thebibliography}
\end{document}